\documentclass{amsart}
\usepackage{amsmath}
\usepackage{amssymb}
\usepackage{amsthm}
\usepackage{amsfonts}
\usepackage{graphicx}
\newtheorem{Theorem}{Theorem}[section]

\newtheorem{Definition}[Theorem]{Definition}
\newtheorem{Corollary}[Theorem]{Corollary}

\newtheorem{Remark}[Theorem]{Remark}
\theoremstyle{definition}

\numberwithin{equation}{section}
\title{Spacelike Loxodromes on Helicoidal Surfaces in Lorentzian n--Space}

\author[M. Babaarslan]{Murat Babaarslan}
\address{Yozgat Bozok University, Department of Mathematics, 66100, Yozgat, Turkey}
\email{murat.babaarslan@bozok.edu.tr}

\author[B. Bekta\c s]{Burcu Bekta\c s Demirci}
\address{Fatih Sultan Vak{\i}f University, Hal\.{I}\c{c} Campus, Faculty of Engineering,
Department of Civil Engineering, 34445, Beyo\u{g}lu, İstanbul, Turkey}
\email{bbektas@fsm.edu.tr}

\author[R. Gen\c{c}]{Rukiye Gen\c{c}}
\address{Yozgat Bozok University, Graduate School of Natural and Applied Sciences, Department of Mathematics, 66100, Yozgat, Turkey}
\email{rkye.gnc@gmail.com}

\begin{document}
\maketitle
\begin{abstract}
In this paper, we study three types of helicoidal surfaces
in a Lorentzian n--space $\mathbb{E}^n_1$. 
First, we find the parametrizations of spacelike loxodromes on such spacelike helicoidal surfaces in $\mathbb{E}^n_1$.
Then, we make a similar classification for spacelike loxodromes on such timelike helicoidal surfaces of a Lorentzian n--space $\mathbb{E}^n_1$. 
\end{abstract}

\section{Introduction}
A loxodrome is a special curve that makes a constant angle with all meridians 
on the Earth's surface. 
Even though the loxodromes do not help finding the shortest
route across the Earth's surface between two points, 
they assures an efficient routing from
one point to another by means of a constant course angle.
Thus, the loxodromes are a popular topic in navigation, see \cite{Alex} and \cite{C}.
From mathematical point of view,
geometers have studied loxodromes defined on different kind of surfaces 
in various ambient spaces 
and then many works have been found in literature such as \cite{BDJ,KFH,COP,BG, BK2, B1, B2, BS, MTV}.

In \cite{Noble}, C. A. Noble obtained the equations of loxodromes on the 
rotational surfaces and he also found such curves on spheres and sphereoids.
In addition, S. Kos et al. \cite{KVZ} calculated the arc--length of the loxodromes on a sphere 
and 
M. Petrovi{\'c} \cite{Petro} also found the arc--length of the loxodromes on a sphereoid.
Moreover, 
there are some obtained results about spacelike loxodromes and timelike
loxodromes on the rotational surfaces in Minkowski 3--space given in \cite{BY1} and \cite{BM}, respectively.
In \cite{Yoon}, D. W. Yoon got results about loxodromes on rotational surfaces 
in the 3--dimensional simply isotropic space. 

Since the helicoidal surfaces are invariant under a motion which is a composition of a translation with a rotation, the helicoidal surfaces are considered a kind of generalization of rotational surfaces.
In \cite{BY2}, M. Babaarslan and Y. Yaylı studied the differential equation of loxodromes on helicoidal surfaces in the 3--dimensional Euclidean space $\mathbb{E}^3$. 
Recently, M. Babaarslan and M. Kayacik studied the spacelike loxodromes on helicoidal surfaces in Minkowski 3--space having spacelike meridians and timelike meridians in \cite{BK} and M. Babaarslan and N. Sönmez \cite{BS} 
also obtained the parametrization of loxodromes on 
the non--degenerate helicoidal surfaces in 
a 4--dimensional Minkowski space
$\mathbb{E}^4_1$. 

In 4--Minkowski space $\mathbb{E}^4_1$, 
there are three types rotation such as spacelike, hyperbolic and screw rotation which leave spacelike,
timelike and lightlike plane respectively.
With this motivation, three types of helicoidal surfaces can be defined in 
n--dimensional Lorentzian space.
In this article, we consider three types of non--degenerate helicoidal surfaces 
in a Lorentzian n--space $\mathbb{E}^n_1$. 
First, we obtain the differential equations of spacelike loxodromes on such spacelike helicoidal surfaces in $\mathbb{E}^n_1$. 
Then, we give result about the parametrizations of
such loxodromes. 
Finally, we make similar computation for spacelike loxodromes 
on timelike helicoidal surfaces of a Lorentzian n--space $\mathbb{E}^n_1$. 
As particular cases, we study spacelike loxodromes 
on non--degenerate right helicoidal surfaces in a Lorentzian n--space $\mathbb{E}^n_1$.

\section{Preliminaries}
For the vectors $x=(x_1, x_2,\dots, x_n)$ and $y=(y_1, y_2, \dots, y_n)$ in $n$--dimensional Euclidean space 
$\mathbb{R}^n$, 
the Lorentzian inner product $\langle\cdot,\cdot\rangle$ 
of $x$ and $y$
is defined by 
\begin{equation}
    \langle x, y\rangle=x_1y_1+x_2y_2+\cdots +x_{n-1}y_{n-1}-x_ny_n.
\end{equation}
$(\mathbb{R}^n, \langle\cdot,\cdot\rangle)$ 
with the Lorentzian inner product given above 
is said to be Lorentzian $n$--space. 
Then, it is denoted by $\mathbb{E}^n_1$. 

The length of the vector $x$ in $\mathbb{E}^n_1$ 
is defined by 
$||x||=\sqrt{|\langle x, x\rangle|}$ 
and it is said to be a unit vector 
if $||x||=1$. 
A causal character of any arbitrary vector $x$ in $\mathbb{E}^n_1$ is said to 
be spacelike(resp., timelike or lightlike) 
if $||x||>0$ or $x=0$
(resp., $||x||<0$ or $||x||=0$ and $x\neq 0$).

Let $\alpha:I\subset\mathbb{R}\longrightarrow\mathbb{E}^n_1$ be a smooth curve
in $\mathbb{E}^n_1$. 
The curve $\alpha$ is spacelike(resp., timelike or lightlike) provided that
$\alpha'$ is spacelike (resp., timelike or lightlike).

Suppose that ${\bf x}:M\longrightarrow\mathbb{E}^n_1$ is an isometric immersion from pseudo--Riemannian surface $M$ to
a Lorentzian n--space $\mathbb{E}^n_1$. 
Then, the coefficients of the first fundamental form of $M$ are 
\begin{equation}
    E=\langle {\bf x}_u, {\bf x}_u\rangle,\;\;
    F=\langle {\bf x}_u, {\bf x}_v\rangle\;\;
    \mbox{and}\;\;
    G=\langle {\bf x}_v, {\bf x}_v\rangle
\end{equation}
for a coordinate system $\{u,v\}$ in $M$. 
Here, ${\bf x}_u$ and ${\bf x}_v$ denote the partial derivative of ${\bf x}$ with respect to $u$ and $v$, 
respectively.
Then, $M$ is a spacelike surface if and only if $EG-F^2>0$; 
$M$ is a timelike surface if and only if $EG-F^2<0$. For $EG-F^2=0$, the surface $M$ is called degenerate. 
Throughout the article, we are not interested in degenerate case. 

Moreover, the length of any non--lightlike curve $\alpha$ on the pseudo--Riemannian surface $M$ in $\mathbb{E}^n_1$
defined by the isometric immersion ${\bf x}$
between two points $u_0$ and $u_1$ is given by 
\begin{equation}
\label{length}
    s=\int_{u_0}^{u_1}\sqrt{\left |E+2F\frac{dv}{du}
    +G\left(\frac{dv}{du}\right)^2\right |}du.
\end{equation}
Assume that $\alpha(t)={\bf x}(u(t),v(t))$ is a spacelike curve on the pseudo--Riemannian surface $M$ in $\mathbb{E}_1^n$, that is, 
\begin{equation}
\label{lox1}
    E\left(\frac{du}{dt}\right)^2+2F\frac{du}{dt}\frac{dv}{dt}
    +G\left(\frac{dv}{dt}\right)^2>0.
\end{equation}
For later use, we also calculate the following equation 
\begin{align}
    \label{lox2}
    \langle\alpha'(t),{\bf x}_u\rangle=E\frac{du}{dt}+F\frac{dv}{dt}.
\end{align}
\begin{Definition}
\cite{Rat}
Let $x$ and $y$ be vectors in a Lorentzian n--space $\mathbb{E}^n_1$. Then, we have the followings:
\begin{itemize}
    \item [i.] for spacelike vectors $x$ and $y$ that span a spacelike vector subspace, 
    there is a unique Lorentzian spacelike angle $\theta$
    between $x$ and $y$ such that 
    \begin{equation}
    \label{ang1}
    \langle x,y\rangle=||x||||y||\cos{\theta},
    \;\;\theta\in[0,\pi],
    \end{equation}
    
    \item [ii.] for spacelike vectors $x$ and $y$ 
    that span a timelike vector subspace, 
    there is a unique Lorentzian timelike angle $\theta$ between $x$ and $y$ such that 
    \begin{equation}
    \label{ang2}
    \langle x,y\rangle=||x||||y||\cosh{\theta},
    \end{equation}
    
    \item [iii.] for a spacelike vector $x$ and a timelike vector $y$ that span a timelike vector subspace,
    there is a unique Lorentzian timelike angle $\theta$
    between $x$ and $y$ such that 
    \begin{equation}
    \label{ang3}
    \langle x,y\rangle=||x||||y||\sinh{\theta}.
    \end{equation}
    \end{itemize}
\end{Definition}
Now, we give the definition of the helicoidal surfaces in $\mathbb{E}^n_1$ as follows.

Assume that 
$\beta:I\subset\mathbb{R}\longrightarrow \Pi\subset\mathbb{E}^n_1$ is a smooth curve in a hyperplane $\Pi\subset\mathbb{E}^n_1$ 
and ${\bf P}$ is $(n-2)$--plane in 
the hyperplane $\Pi\subset\mathbb{E}^n_1$ 
and $\ell$ is a line which does not intersect the curve $\beta$ and is parallel to ${\bf P}$. A helicoidal surface in $\mathbb{E}^n_1$ is defined 
as a rotation of the curve $\beta$ about ${\bf P}$ followed by a translation
along a line $\ell$. Also, the speed of such translation is proportional 
to the speed of this rotation. Hence, we obtain three types of helicoidal surfaces in $\mathbb{E}^n_1$ as follows.

\subsection{Helicoidal surface of type I}
Let $\{e_1,e_2,\dots, e_n\}$ be a standard orthonormal basis for $\mathbb{E}^n_1$,
${\Pi}_I$ a hyperplane spanned by 
$\{e_1, e_3,\dots,e_n\}$ and 
${\bf P}_I$ $(n-2)$--plane spanned by 
$\{e_3, e_4,\dots, e_n\}$. 
Assume that 
$\beta_I:I\longrightarrow\Pi_I\subset\mathbb{E}^n_1,
\;\beta_I(u)=(x_1(u), 0, x_3(u),\dots, x_n(u)),$
is a smooth regular curve in $\Pi_I$ and 
$u$ is arc length parameter, that is, $x_1'^2(u)+x_3'^2(u)+\cdots-x_n'^2(u)=\varepsilon$ with $\varepsilon=\pm 1$.

Then, the parametrization of the helicoidal surface $M_I$ obtained the rotation of the curve $\beta_I$ which leaves the timelike plane ${\bf P}_I$ invariant followed by the translation along $\ell_I$ spanned by $e_n$ can be given 
\begin{align}
\label{typeI}
\begin{split}
&{\bf x}_I(u,v)=
\left[
\begin{array}{ccccc}
\cos{v} & -\sin{v} & 0 & \cdots & 0\\
\sin{v} & \cos{v}  & 0 & \cdots & 0\\
0       &  0       & 1 & \cdots & 0\\
\vdots & \vdots  &  \vdots &  \ddots  & \vdots\\
0 & 0 & 0 & \cdots & 1
\end{array}
\right]
\left[
\begin{array}{c}
     x_1(u)  \\
     0\\
     x_3(u)\\
     \vdots\\
     x_n(u)
\end{array}
\right]
+
cv
\left[
\begin{array}{c}
     0  \\
     0\\
     0\\
     \vdots\\
     1
\end{array}
\right],\\
&\mbox{i.e.,}\\
&M_{I}:{\bf x}_{I}(u,v)=(x_1(u)\cos{v}, x_1(u)\sin{v}, x_3(u), \dots, x_n(u)+cv)
\end{split}
\end{align}
for $0\leq v<2\pi$ and a positive constant $c$. 
By a direct computation, we have
\begin{align}
\begin{split}
\label{typeIbase}
    ({\bf x}_{I})_u&=(x_1'(u)\cos{v}, x_1'(u)\sin{v}, x_3'(u),\cdots, x_n'(u)),\\
    ({\bf x}_{I})_v&=(-x_1(u)\sin{v}, x_1(u)\cos{v}, 0,\cdots, c).
    \end{split}
\end{align}
Thus, the coefficients of the first fundamental form 
of $M_{I}$ is given by
\begin{equation}
    \label{firsttypeI}
    E=\varepsilon,\;\;F=-cx_n'(u)\;\;\mbox{and}\;\; G=x_1^2(u)-c^2.
\end{equation}
Due the fact that the surface $M_{I}$ is non--degenerate 
in $\mathbb{E}^n_1$,
$\varepsilon x_1^2(u)-c^2(\varepsilon+x_n'^2(u))\neq 0$.
If $x_n$ is a constant function, 
then the surface $M_{I}$ is called a 
right helicoidal surface of type I in $\mathbb{E}^n_1$.

\subsection{Helicoidal surface of type II}
Let $\{e_1,e_2,\dots, e_n\}$ be a standard orthonormal basis for $\mathbb{E}^n_1$, 
${\Pi}_{II}$ a hyperplane spanned by 
$\{e_1, e_2,\dots,e_{n-2}, e_n\}$ and 
${\bf P}_{II}$ $(n-2)$--plane spanned by 
$\{e_1, e_2,\dots, e_{n-2}\}$. 
Assume that 
$\beta_{II}:I\longrightarrow\Pi_{II}\subset\mathbb{E}^n_1,
\;\beta_2(u)=(x_1(u), x_2(u),$\\
$\dots,x_{n-2}(u),0, x_n(u)),$
is a smooth regular curve in $\Pi_{II}$
and $u$ is an arc length parameter, that is, 
$x_1'^2(u)+x_2'^2(u)+\cdots-x_n'^2(u)=\varepsilon$ for $\varepsilon=\pm 1$.

Then, the parametrization of the helicoidal surface $M_{II}$ obtained the rotation of the curve $\beta_{II}$ which leaves the spacelike plane ${\bf P}_{II}$ 
invariant followed by the translation along $\ell_{II}$ spanned by $e_1$ can be given 
\begin{align}
\label{typeII}
\begin{split}
&{\bf x}_{II}(u,v)=
\left[
\begin{array}{ccccc}
1 & 0 & \cdots & 0 & 0\\
0 & 1  & \cdots & 0 & 0\\
\vdots &  \vdots & \ddots & \vdots & \vdots\\
0 & 0  &  \cdots &  \cosh{v}  & \sinh{v}\\
0 & 0 & \cdots & \sinh{v} & \cosh{v}
\end{array}
\right]
\left[
\begin{array}{c}
     x_1(u)  \\
     x_2(u)\\
     \vdots\\
     0\\
     x_n(u)
\end{array}
\right]
+
cv
\left[
\begin{array}{c}
     1  \\
     0\\
     0\\
     \vdots\\
     0
\end{array}
\right],\\
&\mbox{i.e.,}\\
&M_{II}:{\bf x}_{II}(u,v)= (x_1(u)+cv, x_2(u), \dots, x_{n-2}(u), x_n(u)\sinh{v}, x_n(u)\cosh{v})
\end{split}
\end{align}
for $v\in\mathbb{R}$ and a positive constant $c$.
By a  direct computation, we get
\begin{align}
\label{typeIIbase}
\begin{split}
    ({\bf x}_{II})_u&=(x_1'(u), x_2'(u), \cdots, x_{n}'(u)\sinh{v}, x_n'(u)\cosh{v}),\\
    ({\bf x}_{II})_v&=(c, 0, \cdots,0, x_n(u)\cosh{v}, x_n(u)\sinh{v}).
    \end{split}
\end{align}
Hence, the coefficients of the first fundamental form 
of $M_{II}$ is given by
\begin{equation}
    \label{firsttypeII}
    E=\varepsilon,\;\;F=cx_1'(u)\;\;\mbox{and}\;\; G=x_n^2(u)+c^2.
\end{equation}
Since $M_{II}$ is a non--degenerate helicoidal surface 
in $\mathbb{E}^n_1$,  
$\varepsilon x_n^2(u)+c^2(\varepsilon-x_1'^2(u))\neq 0$.
In particular, the surface $M_{II}$ is called a right helicoidal surface of type II in $\mathbb{E}^n_1$ when $x_1$ is a constant function.

\subsection{Helicoidal surface of type III}
Let define a pseudo--orthonormal basis 
$\{e_1, e_2, \cdots, \xi_{n-1}, \xi_n\}$ for $\mathbb{E}^n_1$
using a standard orthonormal basis $\{e_1,e_2,\cdots,e_{n-1},e_n\}$ for $\mathbb{E}^n_1$
such that 
\begin{equation}
    \xi_{n-1}=\frac{1}{\sqrt{2}}(e_n-e_{n-1})\;\;\mbox{and}\;\;
    \xi_{n}=\frac{1}{\sqrt{2}}(e_{n}+e_{n-1})
\end{equation}
where $\langle\xi_{n-1},\xi_{n-1}\rangle=\langle\xi_n,\xi_n\rangle=0$ and 
$\langle\xi_{n-1},\xi_n\rangle=-1$. 
Assume that 
$\Pi_{III}$ is a hyperplane spanned by 
$\{e_1, e_3,\cdots,e_{n-2}, \xi_{n-1}, \xi_n \}$ and
${\bf P}_{III}$ is $(n-2)$--plane
spanned by $\{e_1, e_3, \cdots, \xi_{n-1}\}$. 
Suppose that
$\beta_{III}:I\longrightarrow\Pi_{III}\subset\mathbb{E}^n_1,\;\beta_3(u)=x_1(u) e_1+ x_3(u)e_3+\cdots+x_{n-1}(u)\xi_{n-1}+x_n(u)\xi_n$ 
is a smooth curve in $\Pi_{III}$ parameterized by the arc--length parameter $u$, that is, 
$x_1'^2(u)+x_3'^2(u)+\cdots-2x_{n-1}'(u)x_n'(u)=\varepsilon$ 
for $\varepsilon=\pm 1$.

Then, the position vector of 
the helicoidal surface $M_{III}$ obtained a rotation 
of the curve $\beta_{III}$ using the transformation $T$ which leaves the degenerate plane ${\bf P}_{III}$ invariant followed by the translation along $\ell_{III}$ spanned by $\xi_{n-1}$ can be given 
\begin{align}
\label{typeIII}
\begin{split}
 {\bf x}_{III}(u,v)=& x_1(u)e_1+\sqrt{2}vx_n(u)e_2+x_3(u)e_3+\cdots\\
 &+x_{n-2}(u)e_{n-2}
 +(x_{n-1}(u)+v^2x_n(u)+cv)\xi_{n-1}+x_n(u)\xi_n.
 \end{split}
\end{align}
Here, $T$ is an orthogonal transformation 
in $\mathbb{E}^n_1$ such that 
$T(e_1)=e_1,\;T(e_2)=e_2+\sqrt{2}v\xi_{n-1},\; T(e_3)=e_3,\;\dots, 
T(\xi_{n-1})=\xi_{n-1},\;T(\xi_n)=\sqrt{2}ve_2+v^2\xi_{n-1}+\xi_n$. 

By a direct calculation, we find 
\begin{align}
\label{typeIIIbase}
    \begin{split}
       ({\bf x}_{III})_u&=
       x_1'(u)e_1+\sqrt{2}vx_n'(u)e_2+x_3'(u)e_3+\cdots+
       x_{n-2}'(u)e_{n-2}\\
       &+(x_{n-1}'(u)+v^2x_n'(u))\xi_{n-1}
       +x_n'(u)\xi_n,\\
       ({\bf x}_{III})_v&= \sqrt{2}x_n(u)e_2+(2vx_n(u)+c)\xi_{n-1},\\
    \end{split}
\end{align}
and then, the coefficients of the first fundamental form of $M_{III}$ are given by
\begin{equation}
    \label{firstIII}
    E=\varepsilon,\;\;F=-cx_n'(u),\;\;\mbox{and}\;\;
    G=2x_n^2(u).
\end{equation}
Also, $2\varepsilon x_n^2(u)-c^2x_n'^2(u)\neq 0$ because $M_{III}$ is a non--degenerate one.
If $x_n$ is a constant function, then the helicoidal surface $M_{III}$ is called a right helicoidal surface of type III in $\mathbb{E}^n_1$.
\begin{Remark}
It can be easily seen that the helicoidal surfaces $M_{I}$--$M_{III}$
in $\mathbb{E}^n_1$
defined by \eqref{typeI}, \eqref{typeII} and \eqref{typeIII} reduce to
the rotational surfaces in $\mathbb{E}^n_1$ for $c=0$. 
\end{Remark}

\section{Spacelike Loxodromes on Spacelike Helicoidal Surfaces in $\mathbb{E}^n_1$}
In this section, we obtain the parametrization of the spacelike loxodromes on the spacelike helicoidal surface of type I,
type II and type III in a Lorentzian space $\mathbb{E}^n_1$ defined by \eqref{typeI}, \eqref{typeII} and \eqref{typeIII}, respectively.

Let $M_{I}$ be a spacelike helicoidal surface of type I in 
$\mathbb{E}^n_1$ defined by \eqref{typeI} having the spacelike meridians. 
That is, $\varepsilon=1$. 
Then, using the equation \eqref{firsttypeI} the induced metric $g_{I}$ on $M_{I}$ is obtained  
\begin{equation}
    g_{I}=du^2-2cx_n'(u)dudv+(x_1^2(u)-c^2)dv^2
\end{equation}
with $x_1^2(u)-c^2(1+x_n'^2(u))>0$.
Assume that $\alpha_I(t)={\bf x}_{I}(u(t),v(t))$ is 
a spacelike loxodrome on $M_{I}$, 
so that, the equation \eqref{lox1} becomes
\begin{equation}
\label{lox1n}
    \left(\frac{du}{dt}\right)^2
    -2cx_n'(u)\frac{du}{dt}\frac{dv}{dt}
    +(x_1^2(u)-c^2)\left(\frac{dv}{dt}\right)^2>0
\end{equation}
and using the first one of \eqref{typeIbase}, the equation \eqref{lox2} also gives 
\begin{equation}
 \label{lox2n}  
 \langle\alpha'_I(t),({\bf x}_I)_u\rangle=
 \frac{du}{dt}-cx_n'(u)\frac{dv}{dt}.
\end{equation}
Since the loxodrome $\alpha_I(t)$ 
intersects the spacelike meridian 
of $M_I$ with a constant Lorentzian spacelike angle denoting $\theta_0$ at a point $p\in M_{I}$,  
from the equations \eqref{ang1}, \eqref{lox1n} and \eqref{lox2n}, we get
\begin{equation}
  \cos{\theta_0}=\frac{du-cx_n'(u)dv}
  {\sqrt{du^2-2cx_n'(u)dudv+(x_1^2(u)-c^2)dv^2}}.
\end{equation}
After rearranging this equation, we obtain the following differential equation
\begin{equation}
    (\cos^2{\theta_0}(x_1^2(u)-c^2)-c^2x_n'^2(u))
    \left(\frac{dv}{du}\right)^2+2c\sin^2{\theta_0}x_n'(u)
    \frac{dv}{du}=\sin^2{\theta_0}.
\end{equation}
Solving this one, we find the function $v=v(u)$ as  
\begin{equation}
    v(u)=\int_{u_0}^u{\frac{-2c\sin^2{\theta_0}x_n'(\xi)
 \pm\sqrt{(x_1^2(\xi)-c^2(1+x_n'^2(\xi)))
 \sin^2{2\theta_0}}}
 {2\cos^2{\theta_0}(x_1^2(\xi)-c^2)-2c^2x_n'^2(\xi)}}d\xi
\end{equation}
with 
$2\cos^2{\theta_0}(x_1^2(\xi)-c^2)-2c^2x_n'^2(\xi)\neq 0$.
Thus, we get such a parametrization of the spacelike loxodrome 
$\alpha_I(u)={\bf x}_I(u,v(u))$ on the spacelike helicoidal surface $M_{I}$ in $\mathbb{E}^n_1$ given as (i) of Theorem \ref{Thm1}. 

Similarly, 
for a spacelike helicoidal surface of type II $M_{II}$ 
in $\mathbb{E}^n_1$ defined by \eqref{typeII}, we have
$\varepsilon x_n^2(u)+c^2(\varepsilon-x_1'^2(u))>0$ which implies $\varepsilon=1$. 
Then, using the equation \eqref{firsttypeII} the induced metric $g_{II}$ on $M_{II}$ is given by 
\begin{equation}
    g_{II}=du^2+2cx_n'(u)dudv+(x_n^2(u)+c^2)dv^2
\end{equation}
with $x_n^2(u)+c^2(1-x_1'^2(u))>0$.
Assume that $\alpha_{II}(t)={\bf x}_{II}(u(t),v(t))$ 
is a spacelike loxodrome on $M_{II}$.
Hence, the equations \eqref{lox1} and \eqref{lox2} give 
\begin{equation}
\label{lox1n2}
    \left(\frac{du}{dt}\right)^2
    +2cx_n'(u)\frac{du}{dt}\frac{dv}{dt}
    +(x_n^2(u)+c^2)\left(\frac{dv}{dt}\right)^2>0
\end{equation}
and
\begin{equation}
 \label{lox2n2}  
 \langle\alpha'_{II}(t),({\bf x}_{II})_u\rangle=
 \frac{du}{dt}+cx_n'(u)\frac{dv}{dt},
\end{equation}
respectively.
We know that the spacelike curve $\alpha_{II}(t)$ also intersects the spacelike meridian
of $M_{II}$ with a constant Lorentzian spacelike angle at a point $p\in M_{II}$. 
Let say $\theta_0$. 
Using the equations 
\eqref{lox1n2} and \eqref{lox2n2} in \eqref{ang1}, 
we get
\begin{equation}
  \cos{\theta_0}=\frac{du+cx_1'(u)dv}
  {\sqrt{du^2+2cx_1'(u)dudv+(x_n^2(u)+c^2)dv^2}}.
\end{equation}
After doing necessary calculation, we obtain the following differential equation
\begin{equation}
 (\cos^2{\theta_0}(x_n^2(u)+c^2)-c^2x_1'^2(u))
    \left(\frac{dv}{du}\right)^2-2c\sin^2{\theta_0}x_1'(u)
    \frac{dv}{du}=\sin^2{\theta_0}.
\end{equation}
Hence, we find the parametrization of the spacelike loxodrome 
$\alpha_{II}(u)={\bf x}_{II}(u,v(u))$ on the spacelike helicoidal surface $M_{II}$ in $\mathbb{E}^n_1$ 
given as (ii) of Theorem \ref{Thm1}.

Finally, we consider a spacelike helicoidal surface of type III $M_{III}$ in $\mathbb{E}^n_1$ defined by \eqref{typeIII}. 
Then, $2\varepsilon x_n^2(u)-c^2x_n'^2(u)>0$ which says
that $\varepsilon$ must be equal $1$. 
Thus, the equation \eqref{firstIII} gives the induced metric $g_{III}$ on $M_{III}$ as follows 
\begin{equation}
    g_{III}=du^2-2cx_n'(u)dudv+2x_n^2(u)dv^2
\end{equation}
with $2x_n^2(u)-c^2x_n'^2(u)>0$.
Assume that $\alpha_{III}(t)={\bf x}_{III}(u(t),v(t))$ 
is a spacelike loxodrome on $M_{III}$.
Again, the equation \eqref{lox1} implies 
\begin{equation}
\label{lox1n3}
    \left(\frac{du}{dt}\right)^2
    -2cx_n'(u)\frac{du}{dt}\frac{dv}{dt}
    +2x_n^2(u)\left(\frac{dv}{dt}\right)^2>0
\end{equation}
and using the equation \eqref{lox2}, we have 
\begin{equation}
 \label{lox2n3}  
 \langle\alpha'_{III}(t),({\bf x}_{III})_u\rangle=
 \frac{du}{dt}-cx_n'(u)\frac{dv}{dt}.
\end{equation}
Using the the fact that 
and the loxodrome $\alpha_{III}(t)$ 
meets the spacelike meridian
of $M_{III}$ with a constant Lorentzian spacelike angle at a point $p\in M_{III}$, denote $\theta_0$,  
the equations \eqref{ang1}, \eqref{lox1n3} and \eqref{lox2n3} give the following equation
\begin{equation}
  \cos{\theta_0}=\frac{du-cx_n'(u)dv}
  {\sqrt{du^2-2cx_n'(u)dudv+2x_n^2(u)dv^2}},
\end{equation}
which is also expressed as 
\begin{equation}
     (2\cos^2{\theta_0}x_n^2(u)-c^2x_n'^2(u))
    \left(\frac{dv}{du}\right)^2+2c\sin^2{\theta_0}x_n'(u)
    \frac{dv}{du}=\sin^2{\theta_0}.  
\end{equation}
Therefore, we find the parametrization of the  spacelike loxodrome 
$\alpha_{III}(u)={\bf x}_{III}(u,v(u))$ on the spacelike helicoidal surface $M_{III}$ in $\mathbb{E}^n_1$ given as (iii) of Theorem \ref{Thm1}.

Thus, we give the following statement. 
\begin{Theorem}
\label{Thm1}
Let $M$ be a spacelike helicoidal surface in a Lorentzian n--space $\mathbb{E}^n_1$ defined by \eqref{typeI}, \eqref{typeII} and \eqref{typeIII}. Then, the spacelike loxodrome on $M$ has the following parametrization:
\begin{itemize}
    \item [i.]for the helicoidal surface of type I,
    \begin{equation}
     \alpha_I(u)= (x_1(u)\cos{v(u)}, x_1(u)\sin{v(u)}, \dots, x_n(u)+cv(u))
    \end{equation}
    where
    \begin{equation}
 v(u)=\int_{u_0}^u{\frac{-2c\sin^2{\theta_0}x_n'(\xi)
 \pm\sqrt{\sin^2{2\theta_0}(x_1^2(\xi)
 -c^2(1+x_n'^2(\xi)))}}
 {2\cos^2{\theta_0}(x_1^2(\xi)-c^2)-2c^2x_n'^2(\xi)}}d\xi
\end{equation}
provided that 
$\cos^2{\theta_0}(x_1^2(\xi)-c^2)-c^2x_n'^2(\xi)\neq 0$,
    
    \item [ii.] for the helicoidal surface of type II,
    \begin{equation}
    \alpha_{II}(u)= (x_1(u)+cv(u), x_2(u), \dots, x_{n-2}(u), x_n(u)\sinh{v(u)}, x_n(u)\cosh{v(u)})
    \end{equation}
    where
   \begin{equation}
v(u)=\int_{u_0}^u{\frac{2cx_1'^2(\xi)\sin^2{\theta_0}
\pm\sqrt{\sin^2{2\theta_0}(x_n^2(\xi)+c^2(1-x_1'^2(\xi)))}}
{2\cos^2{\theta_0}(c^2+x_n^2(\xi))-2c^2x_1'^2(\xi)}}d\xi
\end{equation}
provided that 
$\cos^2{\theta_0}(c^2+x_n^2(\xi))-c^2x_1'^2(\xi)\neq 0$.

    \item [iii.] for the helicoidal surface of type III,
    \begin{align}
    \begin{split}
     \alpha_{III}(u)=& x_1(u)e_1+\sqrt{2}v(u)x_n(u)e_2+x_3(u)e_3+\cdots+ x_{n-2}(u)e_{n-2}\\
     &+(x_{n-1}(u)+v^2(u)x_n(u)+cv(u))\xi_{n-1}+x_n(u)\xi_n
     \end{split}
    \end{align}
    where
     \begin{equation}
     v(u)=\int_{u_0}^u{\frac{-2cx_n'(\xi)\sin^2{\theta_0}
     \pm\sqrt{\sin^2{2\theta_0}(2x_n^2(\xi)-cx_n'^2(\xi))}}
     {4\cos^2{\theta_0}x_n^2(\xi)-2c^2x_n'^2(\xi)}}d\xi
     \end{equation}
provided that 
$2\cos^2{\theta_0}x_n^2(\xi)-c^2x_n'^2(\xi)\neq 0$ for constants $\theta_0$ and $c>0$.   
\end{itemize}
\end{Theorem}

Using Theorem \ref{Thm1} and the equation \eqref{length}, we get the following Corollary for the spacelike loxodromes on 
the spacelike right helicoidal surfaces in a Lorentzian n--space $\mathbb{E}^n_1$.
\begin{Corollary}
Let $M$ be a spacelike right helicoidal surface in a Lorentzian n--space $\mathbb{E}^n_1$ defined by 
\eqref{typeI}, \eqref{typeII} and \eqref{typeIII}. Then, the length of the loxodrome curve on $M$ between two points $u_0$ and $u_1$ is given by
\begin{equation}
    s=\left|\frac{u_1-u_0}{\cos{\theta_0}}\right|
\end{equation}
where $\theta_0$ is arbitrary constant.
\end{Corollary}

\section{Spacelike Loxodromes on Timelike Helicoidal Surfaces in $\mathbb{E}^n_1$}
In this section, we obtain the parametrization of the spacelike loxodromes on the timelike helicoidal surface of type I,
type II and type III defined by \eqref{typeI}, \eqref{typeII} and \eqref{typeIII}, respectively.

For simplicity, define a constant $\Theta$ as 
\begin{equation}
\label{Theta}
\Theta=\Big\{
\begin{array}{cc}
\cosh{\theta_0} & \mbox{if}\;\;\varepsilon=1,\\
\sinh{\theta_0} & \mbox{if}\;\;\varepsilon=-1.
\end{array}
\end{equation}
Consider a timelike helicoidal surface of type I $M_{I}$ 
in $\mathbb{E}^n_1$ defined by \eqref{typeI} which means 
$\varepsilon x_1^2(u)-c^2(\varepsilon+x_n'^2(u))<0$. 
Thus, there are two cases occur according to the 
casual character of the meridian curve, i.e, $\varepsilon=1$ or $\varepsilon=-1$. 

Assume that $\alpha_I(t)={\bf x}_{I}(u(t),v(t))$ is 
a spacelike loxodrome on $M_{I}$, 
so that, the equation \eqref{lox1} becomes
\begin{equation}
\label{lox1nt}
    \varepsilon\left(\frac{du}{dt}\right)^2
    -2cx_n'(u)\frac{du}{dt}\frac{dv}{dt}
    +(x_1^2(u)-c^2)\left(\frac{dv}{dt}\right)^2>0
\end{equation}
and using the first one of \eqref{typeIbase}, the equation \eqref{lox2} also gives 
\begin{equation}
 \label{lox2nt}  
 \langle\alpha'_I(t),({\bf x}_I)_u\rangle=
 \varepsilon\frac{du}{dt}-cx_n'(u)\frac{dv}{dt}.
\end{equation}
Since the loxodrome $\alpha_I(t)$ 
intersects the meridian 
of $M_I$ with a constant Lorentzian timelike angle
at a point $p\in M_{I}$,  
from the equations \eqref{ang2}, \eqref{ang3},  \eqref{lox1nt} and \eqref{lox2nt}, we get
\begin{equation}
\Theta=
\frac{\varepsilon du-cx_n'(u)dv}
{\sqrt{\varepsilon du^2-2cx_n'(u)dudv+(x_1^2(u)-c^2)dv^2}}.
\end{equation}
After rearranging this equation, we get the following differential equation
\begin{equation}
(\Theta^2(x_1^2(u)-c^2)-c^2x_n'^2(u))
\left(\frac{dv}{du}\right)^2
+2c(\varepsilon-\Theta^2)\frac{dv}{du}
=1-\varepsilon\Theta^2.
\end{equation}
Solving this one, we obtain 
    \begin{equation}
    v(u)=\int_{u_0}^u{\frac{2cx_n'(\xi)
    (\Theta^2-\varepsilon)
    \pm
    \sqrt{\sinh^2{2\theta_0}(c^2(\varepsilon+x_n'^2(\xi))
    -\varepsilon x_1^2(\xi))}}
    {2\Theta^2(x_1^2(\xi)-c^2)-2c^2x_n'^2(\xi)}}d\xi
    \end{equation}
with 
$2\Theta^2(x_1^2(\xi)-c^2)-2c^2x_n'^2(\xi)\neq 0$.
Thus, we find the parametrization of the spacelike loxodrome 
$\alpha_I(u)={\bf x}_I(u,v(u))$ on the timelike helicoidal surface $M_{I}$ in $\mathbb{E}^n_1$ given as (i) of Theorem \ref{Thm2}. 

Similarly, let $M_{II}$ be a timelike helicoidal surface 
of type II in $\mathbb{E}^n_1$ defined by \eqref{typeII}. 
That is 
$\varepsilon x_n^2(u)+c^2(\varepsilon-x_1'^2(u))<0$.
Thus, there are two cases occur according to the 
casual character of the meridian curve, i.e, $\varepsilon=1$ or $\varepsilon=-1$. 
Assume that $\alpha_{II}(t)={\bf x}_{II}(u(t),v(t))$ 
is a spacelike loxodrome on $M_{II}$, 
so that the equation \eqref{lox1}
and \eqref{lox2}
imply 
\begin{equation}
\label{lox1nt2}
    \varepsilon\left(\frac{du}{dt}\right)^2
    +2cx_1'(u)\frac{du}{dt}\frac{dv}{dt}
    +(x_n^2(u)+c^2)\left(\frac{dv}{dt}\right)^2>0
\end{equation}
and 
\begin{equation}
 \label{lox2nt2}  
 \langle\alpha'_{II}(t),({\bf x}_{II})_u\rangle=
 \varepsilon\frac{du}{dt}+cx_1'(u)\frac{dv}{dt},
\end{equation}
respectively. 
Due to the fact that the spacelike curve 
$\alpha_{II}(t)$ intersects the meridians 
of the timelike surface $M_{II}$ 
with a constant Lorentzian timelike angle at a point $p\in M_{II}$,
using the equations \eqref{lox1nt2} and \eqref{lox2nt2} in \eqref{ang1}, we get
\begin{equation}
  \Theta=\frac{\varepsilon du+cx_1'(u)dv}
  {\sqrt{\varepsilon du^2+2cx_1'(u)dudv+(x_n^2(u)+c^2)dv^2}}.
\end{equation}
After rearranging this equation, we get the following differential equation
\begin{equation}
(\Theta^2(x_n^2(u)+c^2)-c^2x_1'^2(u))
\left(\frac{dv}{du}\right)^2
+2c(\Theta^2-\varepsilon)\frac{dv}{du}
=1-\varepsilon\Theta^2.
\end{equation}
Solving this one, we find the parametrization of the spacelike loxodrome 
$\alpha_{II}(u)={\bf x}_{II}(u,v(u))$ on the timelike helicoidal surface $M_{II}$ in $\mathbb{E}^n_1$ given as (ii) of Theorem \ref{Thm2}.

Let $M_{III}$ be a timelike helicoidal surface of type III in $\mathbb{E}^n_1$ defined by \eqref{typeIII}.
Thus, $2\varepsilon x_n^2(u)-c^2x_n'^2(u)<0$ which implies 
$\varepsilon=1$ or $\varepsilon=-1$. 
Suppose that $\alpha_{III}(t)={\bf x}_{III}(u(t),v(t))$ 
is a spacelike loxodrome on $M_{III}$, 
so that the equation \eqref{lox1} implies 
\begin{equation}
\label{lox1nt3}
    \varepsilon\left(\frac{du}{dt}\right)^2
    -2cx_n'(u)\frac{du}{dt}\frac{dv}{dt}
    +2x_n^2(u)\left(\frac{dv}{dt}\right)^2>0
\end{equation}
and 
the equation \eqref{lox2} also gives 
\begin{equation}
 \label{lox2nt3}  
 \langle\alpha'_{III}(t),({\bf x}_{III})_u\rangle=
 \varepsilon\frac{du}{dt}-cx_n'(u)\frac{dv}{dt}.
\end{equation}
Since the spacelike curve $\alpha_{III}(t)$ also meets 
the meridians 
of $M_{III}$ with a constant Lorentzian timelike angle at a point $p\in M_{III}$, 
the equations \eqref{ang2}, \eqref{ang3}, \eqref{lox1nt3} and \eqref{lox2nt3} in \eqref{ang1} give
\begin{equation}
\Theta=\frac{\varepsilon du-cx_n'(u)dv}
{\sqrt{\varepsilon du^2-2cx_n'(u)dudv+2x_n^2(u)dv^2}}
\end{equation}
which implies
\begin{equation}
(2\Theta^2x_n^2(u)-c^2x_n'^2(u))
\left(\frac{dv}{du}\right)^2
+
2c(\varepsilon-\Theta^2)x_n'(u)\frac{dv}{du}
=
1-\varepsilon\Theta^2.
\end{equation}
Hence, we find the parametrization of the spacelike loxodrome 
$\alpha_{III}(u)={\bf x}_{III}(u,v(u))$ 
on the timelike helicoidal surface $M_{III}$ in $\mathbb{E}^n_1$ given as (iii) of Theorem \ref{Thm2}.

Thus, we have the following Theorem.
\begin{Theorem}
\label{Thm2}
Let $M$ be a timelike helicoidal surface in
a Lorentzian n--space $\mathbb{E}^n_1$ defined by \eqref{typeI}, \eqref{typeII} and \eqref{typeIII}. 
Then, the spacelike loxodrome on $M$ has the following parametrization 
\begin{itemize}
    \item [i.] for the helicoidal surface of type I,
    \begin{equation}
     \alpha_I(u)= (x_1(u)\cos{v(u)}, x_1(u)\sin{v(u)}, \dots, x_n(u)+cv(u))
    \end{equation}
    where
    \begin{equation}
    v(u)=\int_{u_0}^u
    {\frac{2cx_n'(\xi)(\Theta^2-\varepsilon)
    \pm
    \sqrt{\sinh^2{2\theta_0}(c^2(\varepsilon+x_n'^2(\xi))
    -\varepsilon x_1^2(\xi))}}
    {2\Theta^2(x_1^2(\xi)-c^2)-2c^2x_n'^2(\xi)}}d\xi
    \end{equation}
    provided that 
    $\Theta^2(x_1^2(\xi)-c^2)-c^2x_n'^2(\xi)\neq 0$,
    
    \item [ii.] for the helicoidal surface of type II,
    \begin{equation}
    \alpha_{II}(u)= (x_1(u)+cv(u), x_2(u), \dots, x_{n-2}(u), x_n(u)\sinh{v(u)}, x_n(u)\cosh{v(u)})
    \end{equation}
    where
     \begin{equation}
     v(u)=\int_{u_0}^u
     {\frac{-2cx_1'(\xi)(\Theta^2-\varepsilon)
     \pm
     \sqrt{\sinh^2{2\theta_0}(c^2(x_1'^2(\xi)-\varepsilon)
     -\varepsilon x_n^2(\xi))}}
     {2\Theta^2(c^2+x_n^2(\xi))-2c^2x_1'^2(\xi)}}dt
     \end{equation}
     provided that
     $\Theta^2(c^2+x_n^2(\xi))-c^2x_1'^2(\xi)\neq  0$,
     
    \item [iii.] for the helicoidal surface of type III,
    \begin{align}
    \begin{split}
     \alpha_{III}(u)=& x_1(u)e_1+\sqrt{2}v(u)x_n(u)e_2+x_3(u)e_3+\cdots+ x_{n-2}(u)e_{n-2}\\
     &+(x_{n-1}(u)+v^2(u)x_n(u)+cv(u))\xi_{n-1}+x_n(u)\xi_n
     \end{split}
    \end{align}
    where
     \begin{equation}
     v(u)=\int_{u_0}^u
     {\frac{2cx_n'(\xi)(\Theta^2-\varepsilon)
     \pm\sqrt{\sinh^2{2\theta_0}(c^2x_n'^2(\xi) 
     -2\varepsilon x_n^2(\xi))}}{4\Theta^2x_n^2(\xi)-2c^2x_n'^2(\xi)}}
     d\xi
     \end{equation}
     
\end{itemize}
provided that 
     $2\Theta^2x_n^2(\xi)-c^2x_n'^2(\xi)\neq 0$ for constants $c>0$ and
$\Theta$ defined by \eqref{Theta}.
\end{Theorem}
Using Theorem \ref{Thm2} and the equation \eqref{length}, we can give the following Corollary for the spacelike loxodromes on the timelike right helicoidal surfaces in a Lorentzian n--space $\mathbb{E}^n_1$.
\begin{Corollary}
Let $M$ be a timelike right helicoidal surface in a Lorentzian n--space $\mathbb{E}^n_1$.  
Then, the length $s$ of a spacelike loxodrome on $M$ 
between two points $u_0$ and $u_1$ is given by the following:
\begin{itemize}
    \item [i.] for a timelike right helicoidal surface of type I, $s=\left|\frac{u_1-u_0}{\Theta}\right|$,
    
    \item [ii.] for timelike right helicoidal surface of type II and III, $s=\left|\frac{u_1-u_0}{\sinh{\theta_0}}\right|$
\end{itemize}
where 
$\Theta$ is defined by \eqref{Theta} and $\theta_0$ is a constant.
\end{Corollary}

\end{document}